\documentclass[12pt, reqno, twoside]{amsart}
\usepackage{amsmath, amsthm, amscd, amsfonts ,, amssymb, graphicx, color}
\usepackage[bookmarksnumbered, colorlinks, plainpages]{hyperref}
\hypersetup{colorlinks=true,linkcolor=red, anchorcolor=green, citecolor=cyan, urlcolor=red, filecolor=magenta, pdftoolbar=true}
\usepackage{enumitem}
\newtheorem{theorem}{Theorem}[section]
\newtheorem{lemma}[theorem]{Lemma}
\newtheorem{corollary}[theorem]{Corollary}
\theoremstyle{definition}
\newtheorem{definition}[theorem]{Definition}

\theoremstyle{remark}
\newtheorem{remark}[theorem]{Remark}
\numberwithin{equation}{section}
\newtheorem{proposition}[theorem]{Proposition}
\begin{document}

\setcounter{page}{1}

\title[Estimates on the first eigenvalue of a quasilinear elliptic system]{Comparison estimates on the first eigenvalue of a quasilinear elliptic system}

\author[A. Abolarinwa]{Abimbola Abolarinwa$^*$}

 \author[S. Azami]{Shahroud Azami}
 
\address{Department of Mathematics,
University of Lagos, Akoka, Lagos State,
Nigeria.}
\email{\textcolor[rgb]{0.00,0.00,0.84}{A.Abolarinwa1@gmail.com}}

\address{Department of Mathematics,
  Faculty of Sciences, Imam Khomeini International University, Qazvin, Iran.}
\email{\textcolor[rgb]{0.00,0.00,0.84}{azami@sci.ikiu.ac.ir}}

\subjclass[2010]{Primary 35P15,  Secondary 47J10; 53C21}

\keywords{$p$-Laplacian, eigenvalue problems, Cheng-type estimates, Faber-Krahn inequality, Cheeger constant.}

\date{April 03, 2020}

\begin{abstract}
 We study a system of quasilinear eigenvalue problems with Dirichlet boundary conditions on complete compact Riemannian manifolds. In particular, Cheng comparison estimates and inequality of Faber-Krahn for the  first eigenvalue of a $(p,q)$-Laplacian are recovered. Lastly, we reprove a Cheeger type estimates for $p$-Laplacian, $1<p<\infty$, from where a lower bound estimate  in terms of Cheeger's constant for the  first eigenvalue of a $(p,q)$-Laplacian is built. As a corollary, the first eigenvalue converges to Cheeger's constant as $p,q\to 1,1.$
\end{abstract}
 \maketitle

\section{Preliminaries and main results}
\subsection{Introduction}
 Let $\Omega$ be a bounded domain in an $N$-dimensional Riemannian manifold $(M,g)$. We prove some comparison estimates of Cheng-type, Cheeger type and Faber-Krahn on the principal Dirichlet eigenvalue for the following quasilinear elliptic system 
\begin{equation}\label{eq11}
\left\{ \begin{array}{l}
\displaystyle \Delta_p u + \lambda |u|^{\alpha-1} |v|^{\beta-1} v = 0  \ \ \hspace{1cm}\ \  \text{in} \ \ \Omega \subseteq M,  \\  
\Delta_q v + \lambda |u|^{\alpha-1} |v|^{\beta-1} u = 0 \ \ \hspace{1cm}\ \  \text{in} \ \ \Omega \subseteq M,  \\  
\displaystyle \hspace{2.19cm}  u = 0,  \ v =0   \ \ \hspace{1.3cm} \text{on} \ \ \partial \Omega,
\\ 
\displaystyle  (u,v) \in W^{1,p}_0(\Omega) \times W^{1,q}_0(\Omega),
\end{array} \right.
\end{equation}
where $1< p, q < \infty$, $\alpha, \beta >0$ are real numbers satisfying $ \alpha/p + \beta/q =1$.

The principal eigenvalue of (\ref{eq11}) denoted by $\lambda_{1,p,q}(\Omega)$ is the least  positive real number for which the system has a nontrivial solution $(u,v)$ called an eigenvector in the product Sobolev
space $W^{1,p}_0(\Omega) \times W^{1,q}_0(\Omega)$ with $u\neq 0$ and $v\neq 0$. The differential operator in (\ref{eq11}) is the so called $p$-Laplacian, that is 
$$\Delta_p u = \text{div}(|du|^{p-2} du) \hspace{.5cm}  \text{for} \ \ u \in W^{1,p}_0,$$
where $\text{div}$ and $d$ are respectively the divergence and gradient operators. When $p=2$, $\Delta_p$ is the usual Laplace-Beltrami operator. The $p$-Laplace operator arises in problems from pure Mathematics such as in the theories of quasiregular and quasiconformal mappings as well as in modelling problems of physical phenomena in non-Newtonian fluids, nonlinear elasticity, glaceology, petroleum extraction, porous media flows and reaction-diffusion processes, See for instance \cite{[Lindq]} for detail description of the $p$-Laplace operator.

\subsection{Eigenvalue problem for the $p$-Laplacian}
The nonlinear eigenvalue problem for the $p$-Laplacian is the following
\begin{align}\label{eq12}
\left\{ \begin{array}{l}
\displaystyle \Delta_p f + \lambda |f|^{p-2}f = 0 \ \ \ \ \ \text{in} \ \ \Omega 
\\ 
\displaystyle  \hspace{2.3cm} f = 0 \ \ \ \ \ \text{on} \ \ \partial \Omega
\end{array} \right.
\end{align}
with $p \in [1,\infty)$. In local coordinates system, $p$-Laplacian is written as
$$\Delta_p f =\frac{1}{\sqrt{|g|}}\frac{\partial}{\partial x^i}( \sqrt{|g|}g^{ij}|\nabla f|^{p-2}\frac{\partial}{\partial x^i} f), $$
where $|g| =\det(g_{ij})$ and $g^{ij}=(g_{ij})^{-1}$ is the inverse metric. The first $p$-eigenvalue $\lambda_{1,p}(\Omega)$ of the $p$-Laplacian is the smallest nonzero number $\lambda$ for which the Dirichlet problem (\ref{eq12}) has a nontrivial solution $f \in W^{1,p}_0(\Omega)$, where the Sobolev space $W^{1,p}_0(\Omega)$ is the completion of $C^\infty_0(\Omega)$ with respect to the Sobolev norm 
$$\|f\|_{1,p} = \Big(\int_\Omega (|f|^p +|\nabla f|^p)d\mu\Big)^{\frac{1}{p}}$$
and $d\mu$ is the Riemannian volume  element of $(M,g)$.
The first $p$-eigenvalue can be variationally characterized by 
\begin{align}\label{eq13}
\displaystyle \lambda_{1,p}(\Omega) = \inf_f \Big\{ \frac{\int_\Omega|df|^p d\mu}{\int_\Omega|f|^p d\mu} \ \ \ \  \Big|  \ \ f\neq 0,  \ \ \ f \in W^{1,p}_0(\Omega)  \Big\}
\end{align}
satisfying the following constraint $\int_\Omega |f|^{p-2} f d\mu=0$. The corresponding eigenfunction is the minimizer of (\ref{eq13}) and satisfies the Euler-Langrage equation
\begin{align}\label{eq14}
\displaystyle \int_\Omega|df|^{p-2}\langle df, d\phi\rangle d\mu - \lambda \int_\Omega|f|^{p-2} \langle f, \phi\rangle d\mu = 0
\end{align}
for $ \phi \in C^\infty_0(\Omega)$ in the sense of distribution. Here and in the rest of the paper, $\langle \cdot, \cdot \rangle$ is the inner product with respect to the metric $g$.  We know that (\ref{eq12}) has weak solutions with only partial regularity in general \cite{[FMST],[Lindq],[Tolk]}.

There have been many interesting geometric results on $ \lambda_{1,p}$ in the recent years, see the follwoing references \cite{[Ab1],[KN03],[LMS05],[Mao],[Mat],[Tak]} for examples. In particular the first author, Abolarinwa \cite{[Ab1]}, Takeuchi \cite{[Tak]}, Matei \cite{[Mat]}, Mao \cite{[Mao]} and  Lima et al \cite{[LMS05]} obtained the classical estimates of Faber-Krahn, Cheeger, Mackean and Cheng-type inequalities on $\lambda_{1,p}$. For evolving manifolds see \cite{[Ab2],[Ab3],[Ab4a],[Ab4],[Ab5]}.

\subsection{Eigenvalue problem for a $(p,q)$-Laplacian}
Here we say that $\lambda$ is an eigenvalue of the system (\ref{eq11}) provided $u \in W^{1,p}_0(\Omega)$ and $v \in W^{1,q}_0(\Omega)$ satisfy the system of Euler-Langrage equations
\begin{align}\label{eq15}
\left\{ \begin{array}{l}
\displaystyle \int_\Omega|du|^{p-2}\langle du, d\phi\rangle d\mu - \lambda \int_\Omega|u|^{\alpha -1} |v|^{\beta -1} \langle v, \phi\rangle d\mu = 0
\\  
\displaystyle  \int_\Omega|dv|^{q-2}\langle dv, d\psi\rangle d\mu - \lambda \int_\Omega|u|^{\alpha -1} |v|^{\beta -1} \langle u, \psi\rangle d\mu = 0
\end{array} \right.
\end{align}
for $\phi \in W^{1,p}_0(\Omega)$ and $\psi  \in W^{1,q}_0(\Omega)$. The  pair $(u,v),\ u>0, v>0$ is the corresponding eigenfunctions. In a similar manner to the first $p$-eigenvalue, the principal $(p,q)$-eigenvalue is variationally characterized as
\begin{align*}
\lambda_{1,p,q}(\Omega) = \inf \{\mathbb{A}(u,v) \ \ \ | \ \ (u,v) \in W^{1,p}_0(\Omega) \times W^{1,q}_0(\Omega), \ \ \mathbb{B}(u,v)=1 \},
\end{align*}
where
$$\mathbb{A}(u,v) = \frac{\alpha}{p} \int_\Omega |du|^p d\mu +  \frac{\beta}{q} \int_\Omega |dv|^q d\mu \ \ \ \text{and} \ \ \mathbb{B}(u,v) = \int_\Omega |u|^{\alpha-1}|v|^{\beta-1}\langle u, v\rangle d\mu $$
for 
$$ \alpha>0, \ \beta>0 \ \ \text{and} \ \ \frac{\alpha}{p}+\frac{\beta}{q}  =1. $$
The existence, simplicity, stability  and some other properties of $\lambda_{1,p,q}(\Omega)$ have been studied in \cite{[AMZ],[BB12],[BP08],[FMST],[GA87],[KMO],[NP]}, see also the references therein. Indeed,  $\lambda_{1,p,q}(\Omega)$ 
 has been proved to be positive and simple for bounded and unbounded domains in $\mathbb{R}^N.$
 Recently,  the second author \cite{[Az]} applied the approach of symmetrization and co-area formula used by the first author \cite{[Ab1]} to obtain some geometric results of Faber-Krahn and Cheeger inequality for $\lambda_{1,p,q}(\Omega)$. This shows that the classical approaches work well for system (\ref{eq11}) without much difficulty involved.

\subsection{Main results}
The major aim of this paper is to prove Cheng-type comparison estimates \cite{[Cheng1],[Cheng2]}, Faber-Krahn-type inequality and Cheeger-type esitmates \cite{[Che],[Mat]} for $\lambda_{1,p,q}(\Omega)$. Precisely, let $B(x_0,r_0)$ be the open geodesic ball of radius $r_0$  centred at $x_0$ in $M$ and $V_N(k,r_0)$ be a geodesic ball of the same radius $r_0$ in an $N$-dimensional space form $M_k$ of constant sectional curvature $k$. Denote the first eigenvalue of (\ref{eq11}) on  $\overline{B(x_0,r_0)}$ by $\lambda_{1,p,q}(B(x_0,r_0))$ and on $\overline{V_N(k,r_0)}$ by $\lambda_{1,p,q}(V_N(k,r_0))$. Then in Section \ref{sec2} we prove the following theorem.
\begin{theorem}\label{thm11}
Let $M$ be an $N$-dimensional complete Riemannian manifold such that its Ricci curvature $Ric(M) \geq (N-1)k, k \in \mathbb{R}$. Then for any $x_0 \in M$ and $r_0 \in (0, d_M)$, where $d_M$ denotes the diameter of $M$, we have 
\begin{align}\label{eq16}
\lambda_{1,p,q}(B(x_0,r_0)) \leq \lambda_{1,p,q}(V_N(k,r_0)).
\end{align}
Equality holds if and only if $B(x_0,r_0)$ is isometric to $V_N(k,r_0)$
\end{theorem}
In a simple language, Cheng's eigenvalue comparison estimate says that when a domain is large, its first Dirichlet eigenvalue is small and the size of the domain accounts for its curvature. A natural consequence of Theorem \ref{thm11} is the following:
\begin{corollary}\label{cor12}
Let $M$ be an $N$-dimensional compact Riemannian manifold with Ricci curvature $Ric(M)\geq (N-1)k,     \ k \in \mathbb{R}$. Then 
\begin{align}
\lambda_{1,p,q}(M) \leq \lambda_{1,p,q}\Big(V_N(k,\frac{d_M}{2})\Big),
\end{align}
 where $d_M$ denotes the diameter of $M$. 
\end{corollary}
In Section \ref{sec3}, we consider the case where $M$ is  compact with positive Ricci curvature $Ric(M) \geq (N-1)k, k>0$. We first prove the Faber-Krahn-type inequality for any domain $\Omega$ in a complete, simply connected Riemannian manifold of constant sectional curvature.
\begin{theorem}\label{thm13}
Let $\Omega$ be a domain and $B(x_0,R)$ be the geodesic ball of radius $R>0$, both in  a complete, simply
connected Riemannian manifold $M_k$ of constant sectional curvature $k$, such that $Vol(\Omega) = Vol(B(x_0,R))$. Then the following inequality holds.
\begin{align}\label{eq19}
\lambda_{1,p,q}(\Omega) \geq  \lambda_{1,p,q}(B(x_0,R))
\end{align}
The equality in (\ref{eq19}) holds if and only if $\Omega$ is the geodesic ball $B(x_0,R)$.
\end{theorem}
Theorem \ref{thm13} says in particular that the geodesic balls are with smallest $\lambda_{1,p,q}$ among all domains of a given volume. A genralization of the Faber-Krahn inequality for $\Delta_p$ on  a ball in the Euclidean $N$-sphere with radius $1/k^2$ has been established by Matei \cite[Theorem 2.1]{[Mat]} (case $p=2$ is due to Berard and Meyer \cite{[BM82]}). Finally, Matei's result for the
 case of $(p,q)$-Laplacian will be discussed at the end of the proof of Theorem \ref{thm13}. In this case, $\Omega$ will be a domain in a compact Riemannian manifold $M$ with positive Ricci curvature $Ric(M) \geq (N-1)k, k>0.$  

To state the last result, we let $h(\Omega)$ be the Cheeger constant defined by
$$h(\Omega):=  \inf_{\Omega'} \frac{Vol_{N-1}(\partial \Omega')}{Vol_N(\Omega')},$$
where the infimum is taken all over open submanifolds $\Omega'$ with compact closure in $\Omega$ and smooth boundary $\partial \Omega'$. 
Here $Vol_{N-1}(\partial \Omega')$ and $Vol_N(\Omega')$ denote the $(N-1)$-dimensional and $N$-dimensional Riemannian volumes  on $\partial \Omega'$  and $\Omega'$, respectively.
\begin{proposition}\label{prop15}
Let $\Omega$ be a bounded domain with smooth boundary in a complete Riemannan manifold. Then
\begin{align}\label{thm110}
\lambda_{1,p,q}(\Omega) \geq \frac{\alpha}{p} \Big(\frac{h(\Omega)}{p}\Big)^p\|u\|_p^p + \frac{\beta}{q} \Big(\frac{h(\Omega)}{q}\Big)^q\|v\|_q^q,
\end{align}
where the pair $(u,v)$ is the corresponding eigenfunctions to $\lambda_{1,p,q}(\Omega)$ and $\|u\|_r$ is $L^r(\Omega)$-norm $\|u\|_{L^r(\Omega)}=(\int_\Omega |u|^r d\mu(x))^{1/r}$. Moreover, $\lambda_{1,p,q}(\Omega)$ converges to Cheeger's constant as $p \to 1$ and $q \to 1$. 
\end{proposition}

\section{Proofs of Theorem  \ref{thm11} and Corollary \ref{cor12}}\label{sec2}
\subsection{Proof of Theorem \ref{thm11}}
\proof
Let there exist a first pair of eigenfunctions, $(\bar{u},\bar{v}), \bar{u}>0, \bar{v}>0$,  of $(p,q)$-Laplacian on  $\overline{V_N(k,r_0)}$ with Dirichlet boundary with $(\bar{u},\bar{v}) \in W^{1,p}_0(\overline{V_N(k,r_0)}) \times W^{1,q}_0(\overline{V_N(k,r_0)})$. Then $\bar{u}, \bar{v}$ are radial (since $V_N(k,r_0)$ is a ball in a simply connected space form which is  two-points homogeneous). 
 Let $r$ be the distance function on $M$ with respect to the point $x_0$, then $(\bar{u} \circ r,\bar{v} \circ r) \in W^{1,p}_0(B(x_0,r_0)) \times W^{1,q}_0(B(x_0,r_0))$  satisfy the boundary conditions. Therefore by definition
  \begin{align}\label{eq21}
\lambda_{1,p,q}(B(x_0,r_0)) \leq \frac{\alpha}{p}\int_{B(x_0,r_0)} |d (\bar{u} \circ r)|^p d\mu + \frac{\beta}{q}\int_{B(x_0,r_0)} |d (\bar{v} \circ r)|^q d\mu 
  \end{align}
  with $\displaystyle \int_{B(x_0,r_0)} |\bar{u}|^{\alpha-1}|\bar{v}|^{\beta-1}\langle \bar{u}, \bar{v}\rangle d\mu = 1$.

Define a $C^\infty$-map $\theta:(0,r_\xi)\times \mathbb{S}^{N-1} \to M$ by $\theta(t,\xi):=\exp_{x_0}(t,\xi),$
where $\mathbb{S}^{N-1}$ is the $(N-1)$-sphere in $T_{x_0}M$ and $\exp_{x_0}$ is a local diffeomorphism from a neighbourhood of $x_0$ in $M$ and $$r_\xi=r_\xi(x_0):=\sup \{t>0 : \exp_{x_0}(s,\xi)  \ \text{is the unique minimal geodesic from} \  x_0 \}.$$
 Since $M$ is complete,
 $B(x_0,r_0)=\{\exp_{x_0}(t,\xi) : \xi \in \mathbb{S}^{N-1} \ \ \text{and}\ \ t \in [0, a(\xi)]\},$ 
where $a(\xi):=\min\{r_\xi, r_0\}$. 
Then integration over $B(x_0,r_0)$ can be pulled back to the tangent space using geodesic polar coordinates. Hence
\begin{align}\label{eq21b}
\left. \begin{array}{l}
\displaystyle  \int_{B(x_0,r_0)} |d(\bar{u} \circ r) |^p d\mu = \int_{\xi \in \mathbb{S}^{N-1}} d\mathbb{S}^{N-1} \int_0^{a(\xi)} \Big|\frac{d \bar{u}}{dt}\Big|^p \times t^{N-1} \theta(t,\xi)dt,
\\ 
\displaystyle  \int_{B(x_0,r_0)} |d(\bar{v}\circ r) |^q d\mu = \int_{\xi \in \mathbb{S}^{N-1}} d\mathbb{S}^{N-1} \int_0^{a(\xi)} \Big|\frac{d\bar{v}}{dt}\Big|^q \times t^{N-1} \theta(t,\xi)dt,
\end{array} \right.
\end{align}
where $d\mathbb{S}^{N-1}$ is the canonical measure of $\mathbb{S}^{N-1} = \mathbb{S}^{N-1}_{x_0}$ and $\theta(t,\xi) \times t^{N-1}=\sqrt{det (g_{ij})}$  is the volume density induced by $\exp_{x_0}$ and $a(\xi) \leq r_0$ such that $\exp_{x_0}(a(\xi),\xi)$ is the cut point of $x_0$ along the geodesic $t \to \exp_{x_0}(t,\xi)$.

Since we have assumed that $\bar{u}$ is everywhere nonnegative and $\bar{u}(r_0)=0$, we have $\displaystyle \frac{d\bar{u}}{dt}(r_0) <0$ in $r_0$-neighbourhood. We have $ \displaystyle \frac{d \bar{u}}{dt} <0$ in $(0,r)$ (see Proposition \ref{prop21} below). Integrating by parts then yields
\begin{align*}
\displaystyle \int_0^{a(\xi)} \Big(-\frac{d\bar{u}}{dt}\Big)^p  t^{N-1} \theta(t,\xi)dt & =  (-\bar{u}) \Big(-\frac{d\bar{u}}{dt}\Big)^{p-1}t^{N-1} \theta(t,\xi) \Bigg|_0^{a(\xi)} - \int_0^{a(\xi)} \frac{(-\bar{u})}{t^{N-1} \theta(t,\xi)} \\
\displaystyle & \times \frac{d}{dt}\Big(t^{N-1} \theta(t,\xi)\Big(-\frac{d\bar{u}}{dt}\Big)^{p-1} \Big)t^{N-1} \theta(t,\xi)dt.
\end{align*}
By a straightforward computation we have
\begin{align*}
\displaystyle \frac{1}{t^{N-1} \theta(t,\xi)} \  & \frac{d}{dt}\Big(t^{N-1} \theta(t,\xi)\Big(-\frac{d\bar{u}}{dt}\Big)^{p-1} \Big)
 \\
\displaystyle & = -\Big(-\frac{d\bar{u}}{dt}\Big)^{p-2} \Bigg[(p-1) \frac{d^2 \bar{u}}{dt^2} +  \Big( \frac{N-1}{t} + \frac{1}{\theta(t,\xi)}\frac{d\theta(t,\xi)}{dt}\Big)\frac{d\bar{u}}{dt} \Bigg].
\end{align*}
Then using the facts that $ \displaystyle \frac{d\bar{u}}{dt}(0)=0$ and $ \displaystyle \Big(-\bar{u} \Big|\frac{d \bar{u}}{dt}\Big|^{p-1}t^{N-1} \theta(t,\xi)\Big)a(\xi)\leq 0$ we obtain
\begin{align}\label{eq22}
\left. \begin{array}{l}
\displaystyle \int_0^{a(\xi)} \Big|\frac{d\bar{u}}{dt}\Big|^p t^{N-1} \theta(t,\xi)dt 
\\ 
\displaystyle \hspace{1cm} \leq - \int_0^{a(\xi)} \bar{u} \Big|\frac{d\bar{u}}{dt}\Big|^{p-2} \Bigg[(p-1) \frac{d^2 \bar{u}}{dt^2} +  \Big( \frac{N-1}{t} + \frac{1}{\theta(t,\xi)}\frac{d\theta(t,\xi)}{dt}\Big)\frac{d\bar{u}}{dt} \Bigg]t^{N-1} \theta(t,\xi)dt.
\end{array} \right.
\end{align}
Notice that
\begin{align*}
\Delta_p \cdot &= \text{div}(|d\cdot|^{p-2}d\cdot)\\
& = |d\cdot|^{p-2} \Delta \cdot + (p-2)|d\cdot|^{p-2} d\cdot d\cdot.
\end{align*}
Since $\bar{u}$ is radial, writing $\Delta$ in geodesic polar coordinates at $k$ we have
$$\Delta \bar{u} = \frac{d^2\bar{u}}{dt^2} + \Big(\frac{N-1}{t} +\frac{1}{\theta_k^N(t,\xi)}\frac{d\theta_k^N(t,\xi)}{dt}\Big)\frac{d \bar{u}}{dt},$$
where $\theta_k^N(t,\xi)$ is the corresponding volume density on $V_N(k,r_0)$ viewed through the exponential map of $M_k$.
Hence 
\begin{align}\label{e23b}
\left. \begin{array}{l}
\displaystyle \Delta_p \bar{u} = \Big|\frac{d\bar{u}}{dt}\Big|^{p-2} \Bigg[\frac{d^2\bar{u}}{dt^2} + \Big(\frac{N-1}{t} +\frac{1}{\theta_k^N(t,\xi)}\frac{d\theta_k^N(t,\xi)}{dt}\Big)\frac{d\bar{u}}{dt} + (p-2) \frac{d^2 \bar{u}}{dt^2}\Bigg]  
\\ \ \\ 
\displaystyle \hspace{1cm} = \Big|\frac{d\bar{u}}{dt}\Big|^{p-2} \Bigg[(p-1) \frac{d^2 \bar{u}}{dt^2} +  \Big( \frac{N-1}{t} + \frac{1}{\theta_k^N(t,\xi)}\frac{d\theta_N^k(t,\xi)}{dt}\Big)\frac{d\bar{u}}{dt} \Bigg]. 
\end{array} \right.
\end{align}
By the assumption on the Ricci curvature and the classical Bishop's comparison theorem, 
\begin{align}\label{eq24}
\frac{d}{dt}\Big(\frac{\theta(t,\xi)}{\theta_k^Nt,\xi)}\Big) \leq 0
\end{align} 
which implies
$$ \frac{d \bar{u}}{dt}\cdot \frac{1}{\theta(t,\xi)}\frac{d\theta(t,\xi)}{dt} \geq \frac{d \bar{u}}{dt}\cdot \frac{1}{\theta_k^N(t,\xi)}\frac{d\theta_k^N(t,\xi)}{dt}.$$
Hence by (\ref{eq22}) we have
\begin{align}\label{eq25}
\left.  \begin{array}{l}
\displaystyle - \bar{u} \Big|\frac{d\bar{u}}{dt}\Big|^{p-2} \Bigg[(p-1) \frac{d^2 \bar{u}}{dt^2} +  \Big( \frac{N-1}{t} + \frac{1}{\theta(t,\xi)}\frac{d\theta(t,\xi)}{dt}\Big)\frac{d\bar{u}}{dt} \Bigg] \\ \ \\
\displaystyle \hspace{2cm} \leq - \bar{u} \Big|\frac{d\bar{u}}{dt}\Big|^{p-2} \Bigg[(p-1) \frac{d^2 \bar{u}}{dt^2} +  \Big( \frac{N-1}{t} + \frac{1}{\theta_k^N(t,\xi)}\frac{d\theta_k^N(t,\xi)}{dt}\Big)\frac{d\bar{u}}{dt} \Bigg] \\ \ \\
\displaystyle \hspace{2cm} = -\bar{u} \Delta_p \bar{u}
\end{array} \right. 
\end{align}
by using definition (\ref{eq24}). Therefore,  combinning  (\ref{eq22}) and (\ref{eq25}) yields
\begin{align}
\displaystyle \int_0^{a(\xi)} \Big|\frac{d \bar{u}}{dt}\Big|^p  \times t^{N-1} \theta(t,\xi)dt  \leq - \int_0^{a(\xi)} \bar{u} \Delta_p \bar{u} dt \Bigg|_{\overline{V_N(k,r_0)}}
\end{align}
and by (\ref{eq21b}) we have 
\begin{align*}
\displaystyle  \int_{B(x_0,r_0)} |d (\bar{u} \circ r) |^p d\mu & \leq -  \int_{\xi \in \mathbb{S}^{N-1}} d\mathbb{S}^{N-1} \int_0^{a(\xi)}\bar{u} \Delta_p \bar{u} dt\\ 
\displaystyle & = - \int_{V_N(k,r_0)} \bar{u} \Delta_p \bar{u} d\mu  \\ 
\displaystyle & = \int_{V_N(k,r_0)} |d (\bar{u} \circ r)|^p d\mu .
\end{align*}
Similarly,
\begin{align*}
\displaystyle  \int_{B(x_0,r_0)} |d( \bar{v} \circ r) |^q d\mu & \leq \int_{V_N(k,r_0)} |d (\bar{v} \circ r)|^q d\mu .
\end{align*}
Then the required inequality in (\ref{eq16}) follows from (\ref{eq21}), that is,
  \begin{align*}
\displaystyle \lambda_{1,p,q}(B(x_0,r_0)) & \leq \frac{\alpha}{p}\int_{B(x_0,r_0)} |d (\bar{u} \circ r)|^p d\mu + \frac{\beta}{q}\int_{B(x_0,r_0)} |d( \bar{v} \circ r)|^q d\mu\\ 
 \displaystyle  & \leq \frac{\alpha}{p}\int_{V_N(k,r_0)} |d (\bar{u} \circ r)|^p d\mu + \frac{\beta}{q}\int_{V_N(k,r_0)} |d (\bar{v} \circ r)|^q d\mu \\
\displaystyle  &  = \lambda_{1,p,q}(V_N(k,r_0)).
  \end{align*}
In conclusion, the equality 
$$\lambda_{1,p}(B(x_0,r_0)) =  \lambda_{1,p}(V_k(k,r_0))$$ 
 holds when there is equality in (\ref{eq22}) and (\ref{eq25}). 
Then we see  that $a(\xi) \equiv r_0$ for almost all $\xi$ and by continuity for all  $\xi$. Hence $\theta(t,\xi) = \theta_k^N(t,\xi)$ which implies equality in  Bishop's inequality (\ref{eq24}). This then proves that $B(x_0,r_0)$ is isometric to $V_N(k,r_0)$.

\qed

The Cheng comparison result is valid regardless of the cut locus, since the Lebesgue measure of the cut locus is $0$ with respect to the $N$-dimensional Lebesgue  measure of the manifold, which implies that integration over the cut locus vanishes.

\begin{proposition}\label{prop21} 
\cite[Proposition 3.1]{[Mao]}
Let $\phi(s)$ be any solution  of 
\begin{align}\label{e28}
[|\phi'(s)|^{p-2} \theta(s)^{N-1}\phi'(s)]' + \lambda \theta(s)^{N-1} |\phi(s)|^{p-2} \phi(s) =0, \hspace{1cm} 1<p<\infty,
\end{align}
where $\theta(s)>0$ on $(0,a)$. Then $\phi'(s) <0$ on $(0,a)$ whenever $\phi(s) >0$ on $(0,a)$ and  $\lambda > 0$. Here $' = \frac{d}{ds}$.
\end{proposition}
\proof
Integrating (\ref{e28})  from $0$ to $s$ yields
$$|\phi'(s)|^{p-2} \theta(s)^{N-1}\phi'(s) = - \lambda \int_0^s \theta(s)^{N-1} |\phi(s)|^{p-2} \phi(s)ds.  $$
The claim of the proposition follows since  $\theta(s)>0$ on $(0,a)$.
\qed
\begin{remark}
Using (\ref{e23b}), functons $\bar{u}$ and $\bar{v}$ which are radial satisfy
\begin{align*}
\left\{ \begin{array}{l}
\displaystyle |\bar{u}'|^{p-2} \Big[(p-1)\bar{u}''  +  \Big( \frac{N-1}{t} + \frac{d\theta_N^k(t,\xi)/dt}{\theta_k^N(t,\xi)}\Big)\bar{u}' \Big] + \lambda |\bar{u}|^{\alpha-1}|\bar{v}|^{\beta-1}\bar{v}=0 \\ \ \\ 
\displaystyle |\bar{v}'|^{q-2} \Big[(q-1)\bar{v}''  +  \Big( \frac{N-1}{t} + \frac{d\theta_N^k(t,\xi)/dt}{\theta_k^N(t,\xi)}\Big)\bar{v}' \Big] + \lambda |\bar{u}|^{\alpha-1}|\bar{v}|^{\beta-1}\bar{u}= 0. 
 \end{array} \right.
\end{align*}
Notice that each equation in the last system is of the form (\ref{e28}), this can be clearly seen by putting $p=q$. 
\end{remark}

\subsection{Proof of Corollary \ref{cor12}}
We can mimick the steps in the proof of Corollary 1.1 in \cite{[Mat]} (see also \cite[Theorem 2.1]{[Cheng1]}) to establish the proof of Corollary \ref{cor12}.

\proof
Let $x_i$ be a point in $M$, such that $B(x_i, \frac{d_M}{2}),\ i =1, 2,...,m$ are pairwise disjoint. Let $r_i$ be the distance function with respect to $x_i$ and $\varphi_i = \varphi \circ r_i$, $\psi_i = \psi \circ r_i$, where $(\varphi, \psi)$ is the first pair of radial eigenfunctions of $V_N(k, \frac{d_M}{2})$. Then by Theorem \ref{thm11}, we have  
  \begin{align*}
\displaystyle \frac{\alpha}{p}\int_{B(x_i, \frac{d_M}{2})} |d \varphi_i|^p d\mu + \frac{\beta}{q}\int_{B(x_i, \frac{d_M}{2})} |d \psi_i|^q d\mu \leq \lambda_{1,p,q}(V_N(k, \frac{d_M}{2}))
  \end{align*}
with $\displaystyle \int_{B(x_i, \frac{d_M}{2})} |\varphi_i|^{\alpha-1}|\psi_i|^{\beta-1}\langle \varphi_i, \psi_i \rangle d\mu = 1.$
 
  We can extend $\varphi_i$  and $\psi_i$ to be zero outside $B(x_i, \frac{d_M}{2})$, then by elementary linear algebra there exist constants $C_i, i=1,2,...,m$ not all equal zero such that 
\begin{align*}
 \int_M  \Big(\sum_{i=1}^{m}C_i\varphi_i\Big)^{\alpha-1} \Big(\sum_{i=1}^{m}C_i\psi_i\Big)^{\beta-1} d\mu =0.
\end{align*}
Since $B(x_i, \frac{d_M}{2})$ are pairwise disjoint , $\sum_{i=1}^{m}C_i\varphi_i \not\equiv 0$ and $\sum_{i=1}^{m}C_i\psi_i \not\equiv 0.$ Hence
 \begin{align*}
\displaystyle \lambda_{1,p,q}(M) & \leq \frac{\alpha}{p}\int_{M} \Big|d \sum_{i=1}^{m}C_i\varphi_i \Big|^p d\mu + \frac{\beta}{q}\int_{M} \Big|d \sum_{i=1}^{m}C_i\psi_i \Big|^q d\mu  \\
\displaystyle &  = \frac{\alpha}{p}\Bigg(\int_{B(x_1, \frac{d_M}{2})} \Big| C_1 d\varphi_1 \Big|^p d\mu + 
 \cdot \cdot \cdot + \int_{B(x_{m}, \frac{d_M}{2})} \Big| C_{m} d\varphi_{m+1} \Big|^p d\mu \Bigg) \\  
 \displaystyle & \hspace{1cm} + \frac{\beta}{q}\Bigg(\int_{B(x_1, \frac{d_M}{2})} \Big|C_1 d \psi_1 \Big|^q d\mu + 
 \cdot \cdot \cdot + \int_{B(x_{m}, \frac{d_M}{2})} \Big| C_{m} d \psi_{m} \Big|^q d\mu \Bigg)\\ 
\displaystyle & = \frac{\alpha}{p}\int_{B(x_i, \frac{d_M}{2})} \Big|d \sum_{i=1}^{m}C_i\varphi_i \Big|^p d\mu + \frac{\beta}{q}\int_{B(x_i, \frac{d_M}{2})} \Big|d \sum_{i=1}^{m}C_i\psi_i \Big|^q d\mu  \\ 
\displaystyle & \leq \lambda_{1,p,q}(V_N(k, \frac{d_M}{2})).
  \end{align*}
  which completes the proof.
\qed
\begin{remark}
Suppose $M$ has nonnegtaive Ricci curvature $Ric(M)\geq 0$, then the above inequality reads
\begin{align*}
\displaystyle \lambda_{1,p,q}(M)  \leq \lambda_{1,p,q}(V_N(0, \frac{d_M}{2})).
\end{align*}
Thus an upper bound can be explicitly found for $ \lambda_{1,p,q}(M)$ by estimating $\lambda_{1,p,q}(V_N(0, \frac{d_M}{2}))$.
The case $p,q=2$ for $\Delta_p$ is found as 
$$\lambda_1(M) \leq \frac{2m^2 N(N+4)}{(d_M)^2}, m=1$$
by Cheng in \cite[Corollary 2.2]{[Cheng1]}.
\end{remark}

\section{Faber-Krahn type inequality} \label{sec3}
The main tools that will be employed in the proof of Theorem \ref{thm13} are symmetrization procedure and inequalities of P\'olya-Szeg\"o and Hardy-Littlewood. We first recall the definition of the symmetrization of
a function and its properties
\begin{definition}
Let $\Omega$ be an open bounded subset of $\mathbb{R}^N$. Let $f$ be a nonnegative measurable function in $\Omega$ which vanishes on the boundary. The set $\{x \in \Omega : f(x) >t, t>0\}$ is called the level set of $f$.
Let $\Omega^*$ be the ball centred at the origin in $\mathbb{R}^N$ with the same volume as $\Omega$. The function $f^*:\Omega^* \to \mathbb{R}^+$ is the nonincreasing symmetric rearrangement of $f$ with
$$ Vol \{x \in \Omega^* : f^*(x) >t \} = Vol\{x \in \Omega : f(x) >t \}$$ 
\end{definition}

For more details on symmetrization see \cite{[LL],[Tah]}. Now using the above symmetrization we have the following Lemma
\begin{lemma}\label{lem32}
Let $\Omega$ be a compact domain. Let $f, g : \Omega \to \mathbb{R}^+$ be a nonnegative measurable functions and $f^*, g^*: B(x,R) \to \mathbb{R}^+$ be radially nonincreasing functions such that $Vol(\Omega)=Vol(B(x,R))$. Then
\begin{enumerate}
\item $\int_\Omega f = \int_{B(x,R)} f^*$ - equimeasurability of level sets.
\item $\int_\Omega |df|^p \geq \int_{B(x,R)} |df^*|^p$,  \ \ for \ $p>1$ - P\'olya-Szeg\"o inequality. 
\item $\int_\Omega f g \leq  \int_{B(x,R)} f^* g^*$  - Hardy-Littlewood inequality.
\end{enumerate} 
\end{lemma}

P\'olya-Szeg\.o inequality says the Dirichlet integral $\int_\Omega |df|^p$ decreases under the influence of symmetrization. This inequality can be realized by combining the co-area formula and H\"olders inequality. A version of the proof of (1) and (2) is in Aubin \cite[Proposition 2.17]{[Aub]}. The Hardy-Littlewood inequality's proof is contained in \cite[Theorem 2.2, p. 44]{[BS98]}. Though, the nonincreasing rearrangement does not preserve product of functions in general, the equality in (3) is attainable for suitable functions, see \cite{[BS98]} for details.
 
\subsection{Proof of Theorem \ref{thm13}}
Let $u^*$ and $v^*$ be the nonincreasing rearrangement of $u \in W^{1,p}(\Omega)$ and $v \in W^{1,q}(\Omega)$ respectively. Let $(u,v)$ be the minimizing eigenfunction satisfying \\ $\int_\Omega |u|^{\alpha-1} |v|^{\beta-1}\langle u, v\rangle d\mu= 1$ for $\alpha, \beta >0$. We know that $u>0$, \ $v>0$. It follows from P\'olya-Szeg\.o inequality that 
$$\int_\Omega |du|^p d\mu \geq \int_{B(x_0,R)} |du^*|^p d\mu\ \ \ \text{and} \ \ \ \int_\Omega |dv|^q d\mu \geq \int_{B(x_0,R)} |dv^*|^qd\mu.$$
Then
\begin{align*}
\displaystyle \lambda_{1,p,q}(\Omega) & = \frac{\alpha}{p} \int_\Omega |du|^p d\mu + \frac{\beta}{q} \int_\Omega |dv|^q d\mu \\ 
\displaystyle & \geq  \frac{\alpha}{p} \int_{B(x_0,R)} |du^*|^p d\mu + \frac{\beta}{q} \int_{B(x_0,R)} |dv^*|^qd\mu. 
\end{align*}
On the other hand it is clear  from Hardy-Littlewood inequality in Lemma \ref{lem32} that 
$\int_\Omega u^\alpha v^\beta d\mu \leq \int_{B(x_0,R)} u^{*\alpha} v^{*\beta} d\mu$ for real numbers $\alpha, \beta >0$,
where $u^* \in W^{1,p}(B(x_0,R))$ and $v^* \in W^{1,q}(B(x_0,R))$. Note also that one can identify $\int_\Omega u^\alpha v^\beta d\mu$ with $\int_\Omega |u|^{\alpha-1} |v|^{\beta-1}\langle u, v\rangle d\mu$. Clearly since $u>0, v>0$, then 
\begin{align*}
{B}(u,v) & =\int_\Omega |u|^{\alpha-1} |v|^{\beta-1}\langle u, v\rangle d\mu = \int_\Omega u^\alpha v^\beta d\mu \\  & \leq \int_{B(x_0,R)} u^{*\alpha} v^{*\beta} d\mu= \int_{B(x_0,R)} |u^*|^{\alpha-1} |v^*|^{\beta-1}\langle u^*, v^*\rangle d\mu = \mathbb{B}(u^*,v^*).
\end{align*}
We therefore conclude that
\begin{align*}
\displaystyle \frac{\alpha}{p} & \int_{B(x,R)} |du^*|^p d\mu  + \frac{\beta}{q} \int_{B(x,R)} |dv^*|^qd\mu \\
\displaystyle & \geq \inf \{\mathbb{A}(u^*,v^*) \  | \ (u^*,v^*) \in W^{1,p}_0(B(x,R)) \times W^{1,q}_0(B(x,R)), \ \ \mathbb{B}(u^*,v^*)=1 \} \\
\displaystyle & = \lambda_{1,p,q}(B(x,R)).
\end{align*}
Equation (\ref{eq19}) is therefore proved.

\qed

Matei's result \cite{[Mat]} for the case of $(p,q)$-Laplacian can be stated as follows

\begin{theorem}
Let $M$ be an $N$-dimensional compact Riemannian manifold with $Ric(M)\geq (N-1)k$ and $B(R)$ be a geodesic ball of Radius $R>0$ in the Euclidean $N$-sphere $\mathbb{S}^N_k$ with constant sectional curvature $k$ such that for a bounded domain $\Omega \in M$   
$$\frac{Vol(\Omega)}{Vol(M)} =\frac{Vol(B(R))}{Vol(\mathbb{S}^N_k)}.$$
Then
\begin{align}\label{eq31}
\lambda_{1,p,q}(\Omega) \geq \gamma \lambda_{1,p,q}(B(R)),
\end{align}
where $\gamma = Vol(M)/Vol(\mathbb{S}^N_k)$. There is equality if and only if there is an isometry which sends $\Omega$ to $B(R) \subset \mathbb{S}^N_k$.
\end{theorem}
The inequality (\ref{eq31}) reduces to (\ref{eq19}) of Theorem \ref{thm13} if $M =\mathbb{S}^N_k$. 

\section{Lower bound estimates (Cheeger-type)}

In this section, we want to give a lower bound for $\lambda_{1,p,q}(\Omega)$ in terms of the so-called Cheeger's constant and a lower bound estimate on the first eigenvalue of $p$-Laplacian.

\begin{definition}\label{def51}
The Cheeger's  constant $h(\Omega)$ of a domain  $\Omega$ is defined to be
\begin{align}\label{eq51}
h(\Omega):=  \inf_{\Omega'} \frac{Vol_{N-1}(\partial \Omega')}{Vol_N(\Omega')},
\end{align}
where $\Omega'$ ranges  over smooth subdomains of $\Omega$ with compact closure in $\Omega$ with smooth boundary $\partial \Omega'$, and 
 $Vol_{N-1}(\partial \Omega')$ and $Vol_N(\Omega')$  denote the 
 volumes  of $\partial \Omega'$  and $\Omega'$, respectively.
\end{definition}

Let $D$ vary over all smooth subdomain of $\Omega$ whose boundary $\partial D$ does not touch $\partial \Omega$, the quantity $\mathcal{Q}(D):= Vol(\partial D)/Vol(D)$ is called the Cheeger quotient of $D$. Any subdomain $E \subset \Omega$ which realizes the infimum in (\ref{eq51}) is referred to as Cheger domain in $\Omega$ while $\Omega$ is called self-Cheeger if it is a minimizer. Problems involving Cheeger's constant/domains are very interesting in Geometric Analysis. For existence, (non)uniqueness and regularity of Cheeger domains see \cite{[KF03]}. For introductory survey and some physical applications of Cheeger's constant see \cite{[Par]} and the references therein, and see \cite{[Ben]} for further results in manifold setting.

\begin{theorem}\label{thm52}(Cheeger type estimate)
Let $\Omega$ be a bounded domain with smooth boundary in a complete Riemannan manifold. Then
\begin{align}
\lambda_{1,p}(\Omega) \geq \Big(\frac{h(\Omega)}{p}\Big)^p,  \hspace{1cm} 1<p<\infty.
\end{align}
\end{theorem}
The above theorem was originally proved by Cheeger \cite{[Che]} for $p=2$ in the case of manifolds without boundary, and an extension for the general $p$ was given by Matei \cite{[Mat]}. See also \cite{[Ab1],[KF03],[Tak],[Mao]} for the general manifolds and $p>1$.
For completeness and the importance of Theorem \ref{thm52} to the proof of Proposition \ref{prop15}, we repeat the proof here. 

\proof (Theorem  \ref{thm52})
Suppose $ \varphi \in C^\infty_0(\Omega)$ is a positive function and we let $A(t):= \{x \in \Omega : \varphi(x)>t\}$ and $\partial A(t) := \{x \in \Omega : \varphi(x)=t\}$. Using the co-area formula 
\begin{align*}
\displaystyle \int_\Omega |d \varphi| d\mu &=\int_{-\infty}^{\infty}\Big(\int_{A(t)} d A(t)\Big) dt = \int_{-\infty}^{\infty}Vol_{N-1}(\partial A(t)) dt \\ 
\displaystyle & = \int_{-\infty}^{\infty} \frac{Vol_{N-1}(\partial A(t))}{Vol_N(A(t))} \cdot Vol_N(A(t)) dt \\ 
\displaystyle & \geq  \inf_{\Omega' \subset \subset \Omega} \frac{Vol_{N-1}(\partial \Omega')}{Vol_N(\Omega')} \int_{-\infty}^{\infty} Vol_N(A(t)) dt  
 = h(\Omega) \int_\Omega \varphi(x) d\mu.
\end{align*}
The above condition also holds for $\varphi \in W_0^{1,1}(\Omega)$ since $C^\infty_0(\Omega)$ is dense in $W_0^{1,1}(\Omega) $. Now for any $p>1$ and $u \in W_0^{1,p}(\Omega)$, define $\Phi(u):=u^p$. Then, by H\"older's inequality
\begin{align}\label{eq53}
\displaystyle \int_\Omega |d\Phi(u)| d\mu = \int_\Omega |d u^p| d\mu \leq p\Big( \int_\Omega |u|^p d\mu\Big)^{(p-1)/p} \Big( \int_\Omega |d u|^p d\mu\Big)^{1/p}.
\end{align}
Letting $\varphi =u^p$, we have by (\ref{eq53}) 
\begin{align*}
\displaystyle h(\Omega) \leq \frac{\int_\Omega |d \varphi| d\mu}{\int_\Omega |\varphi| d\mu} =  \frac{\int_\Omega |d u^p| d\mu}{\int_\Omega |u^p| d\mu} & \leq \frac{p\Big( \int_\Omega |u|^p d\mu\Big)^{(p-1)/p} \Big( \int_\Omega |d u|^p d\mu\Big)^{1/p}}{\int_\Omega |u|^p d\mu} \\
\displaystyle & = p\Bigg(\frac{\int_\Omega |d u|^p d\mu}{\int_\Omega |u|^p d\mu}\Bigg)^{1/p}.
\end{align*}
Since $u \in W_0^{1,p}(\Omega)$ was arbitrary we arrive at 
\begin{align*}
\frac{h(\Omega)}{p}\leq \lambda_{1,p}^{1/p}(\Omega),
\end{align*}
which concludes the result.

\qed

\begin{corollary} (\cite{[KF03]}
The first eigenvalue of $p$-Laplacian $\lambda_{1,p}(\Omega)$ converges to Cheeger's constant $h(\Omega)$ as $p\to 1$.
\end{corollary}
The proof is in \cite{[KF03]} and we omit it here. This corollary simply implies that if we take $\lambda_{1,1}(\Omega):= \lim \sup_{p\to 1} \lambda_{1,p}(\Omega) = h(\Omega),$ Then one asks for the solvability of the limiting problem
\begin{align}\begin{aligned}
\displaystyle - div \Big(\frac{du}{|du|}\Big)& = \lambda_{1,1}(\Omega) \hspace{1cm}  \text{in} \ \ \ \ \Omega \\
\displaystyle u & = 0 \hspace{2cm} \text{on} \ \ \ \ \partial \Omega.
\end{aligned} \end{align}

The proof of Proposition \ref{prop15} is based on the proof of Theorem \ref{thm52}, it is therefore summarised below.

\subsection*{Proof of Proposition \ref{prop15}}
Let $(u,v) \in  W_0^{1,p}(\Omega) \times \in W_0^{1,q}(\Omega)$ be the pair of eigenfunctions corresponding to $\lambda_{1,p,q}(\Omega)$ with $u>0, \ v>0$ by defintion. Then 
$$\lambda_{1,p,q}(\Omega) =  \frac{\alpha}{p} \int_\Omega |du|^p d\mu +  \frac{\beta}{q} \int_\Omega |dv|^q d\mu \ \
\text{with} \ \ 
\int_\Omega |u|^{\alpha-1}|v|^{\beta-1}\langle u, v\rangle d\mu =1. $$
Using Theorem \ref{thm52} we have
\begin{align}
\lambda_{1,p}(\Omega) = \frac{\int_\Omega |d u|^p d\mu}{\int_\Omega |u|^p d\mu} \geq \Big(\frac{h(\Omega)}{p}\Big)^p
\end{align} 
which implies
$$\int_\Omega |d u|^p d\mu \geq  \Big(\frac{h(\Omega)}{p}\Big)^p \|u\|_p^p  \ \ \ \ \ \text{and} \ \ \ \ \ \ \int_\Omega |d v|^q d\mu \geq  \Big(\frac{h(\Omega)}{p}\Big)^p \|v\|_q^q $$
since $u$ and $p$ in Theorem \ref{thm52} were abitrary and for the compact embedding of $ W_0^{1,p}(\Omega) \hookrightarrow  L^p(\Omega)$. Therefore
\begin{align*}
\displaystyle \lambda_{1,p,q}(\Omega)=  \frac{\alpha}{p} \int_\Omega |du|^p d\mu +  \frac{\beta}{q} \int_\Omega |dv|^q d\mu \geq \frac{\alpha}{p} \Big(\frac{h(\Omega)}{p}\Big)^p\|u\|_p^p + \frac{\beta}{q} \Big(\frac{h(\Omega)}{q}\Big)^q\|v\|_q^q, 
\end{align*}
which proves (\ref{thm110}).
Next thing to do is to obtain the limiting behaviour of  $\lambda_{1,p,q}$ as $p \to 1$ and $q \to 1$. Heuristically as $p \to 1$ and $q \to 1$, $\alpha/p + \beta/q \to \alpha+\beta=1$. Taking $\|u\|_r =1$, i.e., we normalize the eigenfunction of $r$-Laplacian and then observe that the lower bound in the last inequality converges to $h(\Omega)$:
$$ \limsup_{p\to 1, \\ q\to 1}\lambda_{1,p,q}(\Omega) = h(\Omega).$$
It therefore suffices to obtain a finite bound for  $\lambda_{1,p,q}(\Omega)$ as $p \to 1$ and $q \to 1$. To obtain a suitable upper bound for 
$\lambda_{1,r}(\Omega)$ (resp. $\int_\Omega |du|^r d\mu, \ r>1$), we can follow the proof of Corollary $6$ of \cite{[KF03]} and then conclude that 
$$\lambda_{1,1,1}(\Omega):= \limsup_{p\to 1, \\ q \to 1} \lambda_{1,p,q}(\Omega) = h(\Omega).$$


\begin{thebibliography}{9}
 \bibitem {[Ab1]}
			{A. Abolarinwa},  
			{The first eigenvalue of $p$-Laplacian and geometric estimates}, {\it Nonl. Anal. Diff. Eq.}, 2(3)(2014), 105--115

\bibitem {[Ab2]}
			{A. Abolarinwa}, 
			{Evolution and monotonicity of the first eigenvalue of $p$-Laplacian under the Ricci-harmonic flow}, {\it J. Appl. Anal.}, 21(2)(2015), 147--60.
			

\bibitem {[Ab3]}
			{A. Abolarinwa, O. Adebimpe and E. A. Bakare}, 
			{Monotonicity formulas for the first eigenvalue of the weighted $p$-Laplacian under the Ricci-harmonic flow},  {\it J. Ineq. Appl.}, (2019) 2019:10. 
			
 \bibitem {[Ab4a]}
			{A. Abolarinwa}, 
			{Eigenvalues of weighted-Laplacian under the extended Ricci flow}, {\it Adv. Geom.}, 19(1) (2019), 131--143.  
			
\bibitem {[Ab4]}
			{A. Abolarinwa, S. O. Edeki and J. Ehigie}, 
			{On the spectrum of the weighted $p$-Laplacian under the Ricci-harmonic flow},  {\it J. Ineq. Appl.}, (2020)(58), 1--14. 			
\bibitem {[Ab5]}
			{A. Abolarinwa, C. Yang and D. Zhang}, 
			{On the spectrum of the $p$-biharmonic operator under the Ricci flow},  {\it Res. Math.}, (2020)(54), 1--16. 
			
			
\bibitem {[Aub]}
			{T. Aubin}, 
			\emph{Some nonlinear problems in Riemannian geometry}, Springer, Berlin (1998).
			
			
\bibitem {[AMZ]}
			{G. A. Afrouzi, M. Mirzapour and Q. Zhang}, {Simplicity and stablity of the first eigenvalue
of a $(p,q)$ Laplacian  system}, {\it Electr. J. Diff Eqn.},  2012(8)(2012), 1--6.


\bibitem {[Az]}
			{S. Azami}, 
			{The first eigenvalue of some $(p, q)$-Laplacian and geometric estimates}, {\it Commun. Korean Math. Soc.},  33(1) (2018), 317--323.
			
\bibitem {[BS98]}
			{C. Bennett and R. Sharpley}, 
			\emph{Interpolation of Operators},
Vol. 120 Pure and Applied Math., Academic press Inc. (1998).
			
\bibitem {[BB12]}
			{N. Benouhiba and Z. Belyacine}, 
			{A class of eigenvalue problems for the $(p, q)$-Laplacian in $\mathbb{R}^N$},
{\it Int. J. of pure appl. Math.},  80(5)(2012), 727--737.


\bibitem {[Ben]}
			{B. Benson}, 
			{The Cheeger constant, isoperimetric problems, and hyperbolic surfaces}, arxiv.org/abs/1509.08993

\bibitem{[BM82]}
			{P. Berard and D. Meyer}, 
			{Inegalities isoperimetriques et applications}, {\it Ann. Scient. Ec. Norm. Sup. 4e Ser.}
15 (1982) 513--542.


\bibitem{[BP08]}
	 {F. Bonder and J. P.  Pinasco},
	 {Estimates for eigenvalues of quasilinear elliptic systems. Part II}, {\it J. Diff. Eq.} 245 (2008), 875--891

			
\bibitem {[Che]}
			{J. Cheeger}, 
			{A lower bound for the smallest eigenvalue of the Laplacian. In Problems in analysis}, (Papers dedicated to Salomon Bochner, 1969), pp. 195--7199. Princeton Univ. Press, Princeton, N. J., 1970.
			
	
\bibitem{[Cheng1]}
	{S. Y. Cheng},
	{Eigenvalue comparison theorems and its geometric applications}, {\it Math. Z.}  141(1975), 289--297.

	
\bibitem{[Cheng2]}
	{S. Y. Cheng},
	{Eigenfunctions and eigenvalues of Laplacian}, {\it Proc. Symp. Pure Math.}, 27(1975), 185--193.

\bibitem{[FMST]} 
		{J. Fleckinger, R. F. Man\'asevich, N. M. Stavrakakis and F. de Th\'elin}, 
		{Principal eigenvalues for some quasilinear elliptic equations on $\mathbb{R}^N$},  {\it Adv. Diff. Eq.} 2(6) (1997), 981--1003.

	
\bibitem{[GA87]}
	{J. P. Garzia Azorero and I. Peral Alonso}, 
	{Existence and nonuniqueness for the p-Laplacian eigenvalues}, {\it Comm. Part. Diff. Eq.} 12(1987) 1389--1430.
	
	     
			       		
\bibitem {[KN03]}
			{S. Kawai and N. Nakauchi}, 
			{The first eigenvalue of the $p$-Laplacian on 		a compact manifold}, {\it Nonl. Anal.}, 55 (2003), 33--46.
			
\bibitem {[KF03]}
			{B. Kawohl and V. Fridman}, 
			{Isoperimetric estimates for the first 			eigenvalue of the p-Laplace operator and the Cheeger constant}, {\it Comment. Math. Univ. Carolin.}, 44 (2003) 659--667.			
			
\bibitem {[KMO]}
			{A. E. Khalil, S. E. Manouni, M. Ouanan}, 
			{Simplicity and stablity of the first eigenvalue
of a nonlinear elliptic system}, {\it Int. J. Math. and Math. Sci.}, 10 (2005) 1555--1563.
			
	
\bibitem {[LL]}
			{E. H. Lieb and M. Loss}, 
			\emph{Analysis}, Graduate Studies in Mathematic, Vol. 14  AMS, Providence, Rhode Island (2007).
			
			
\bibitem {[LMS05]}
			{B. P. Lima, J. F. Montenegro and N. L. Santos}, 
			{Eigenvalues estimates for the first eigenvalue of the $p$-Laplace operator on manifolds}, arXiv:0808.2028[math.DG] 14 Aug 2008.

\bibitem{[Lindq]}
{P. Lindqvist},
{On the equation $div(|\nabla u|^{p-2}\nabla u) + |u|^{p-2} u = 0$}, {\it Res. Reports}c, vol. A 263, Helsinki
Univ. Tech. Inst. Math., Helsinki, 1988.


\bibitem {[Mao]}
			{J. Mao},
			{Eigenvalue inequalities for the p-Laplacian on a Riemannian manifold and estimates for the heat kernel}, {\it J. Math. Pures Appl.} 101 (2014) 372–393

\bibitem {[Mat]}
			{A. M. Matei}, 
			{First eigenvalue for the p -Laplace operator}, {\it Nonl. Anal.} 39(2000) 1051--1068.
			
 
\bibitem {[NP]}
			{P. L. De N\'apoli and J. P. Pinasco}, 
			{Estimates for eigenvalues of quasilinear elliptic systems}, {\it J. Diff. Eq.} 227 (2006), 102--115.
			
\bibitem {[Par]}
			{E. Parini}, 
			{An introduction to the Cheeger problem}, {\it Surveys Math. Appl.}, 6 (2011), 9--22. 
			

\bibitem{[Tah]}
		{ A. Taheri},
		\emph{Function spaces and partial differential equations I \& II}, Oxford Lecture Series in Mathematics and its Applications 40 \& 41,  Oxford University Press, 2015.
		
											
\bibitem {[Tak]}
			{H. Takeuchi}, 
			{On the first eigenvalue of the $p$-Laplacian in a Riemannian manifold}, {\it Tokyo J. Math.}, 21(1)(1998), 136--140.

	\bibitem{[Tolk]}
	{P. Tolksdorf},
	{Regularity for a more general class of quasilinear elliptic equations}, {\it Jour. of Diff. Eq.}, 5(1)(1984), 126--150.
\end{thebibliography}
\end{document}